\documentclass[11pt]{amsart}
\usepackage[all]{xy}
\begin{document}
\title{Some boundedness results for Fano-like Moishezon manifolds}
\author{Laurent BONAVERO and Shigeharu TAKAYAMA}
\date{January 2000}
\maketitle
\noindent
\def\restriction{\string |}
\newcommand{\pp}{\rm ppcm}
\newcommand{\pg}{\rm pgcd}
\newcommand{\Ker}{\rm Ker}
\newcommand{\C}{{\mathbb C}}
\newcommand{\Q}{{\mathbb Q}}
\newcommand{\GL}{\rm GL}
\newcommand{\SL}{\rm SL}
\newcommand{\diag}{\rm diag}
\newcommand{\N}{{\mathbb N}}
\def\finpreuve
{\hskip 3pt \vrule height6pt width6pt depth 0pt}

\newtheorem{theo}{Theorem}
\newtheorem{prop}{Proposition}
\newtheorem{lemm}{Lemma}
\newtheorem{lemmf}{Lemme fondamental}
\newtheorem{defi}{Definition}
\newtheorem{exo}{Exercice}
\newtheorem{rem}{Remark}
\newtheorem{cor}{Corollaire}
\newcommand{\CC}{{\mathbb C}}
\newcommand{\ZZ}{{\mathbb Z}}
\newcommand{\RR}{{\mathbb R}}
\newcommand{\QQ}{{\mathbb Q}}
\newcommand{\FF}{{\mathbb F}}
\newcommand{\PP}{{\mathbb P}}
\newcommand{\NN}{{\mathbb N}}
\newcommand{\codim}{\operatorname{codim}}
\newcommand{\Ho}{\operatorname{Hom}}
\newcommand{\Pic}{\operatorname{Pic}}
\newcommand{\NE}{\operatorname{NE}}
\newcommand{\Nun}{\operatorname{N}}
\newcommand{\card}{\operatorname{card}}
\newcommand{\Hilb}{\operatorname{Hilb}}
\newcommand{\mult}{\operatorname{mult}}
\newcommand{\vol}{\operatorname{vol}}
\newcommand{\divi}{\operatorname{div}}
\newcommand{\pr}{\operatorname{pr}}
\newcommand{\con}{\operatorname{cont}}
\newcommand{\ima}{\operatorname{Im}}
\newcommand{\rk}{\operatorname{rk}}
\newcommand{\Exc}{\operatorname{Exc}}
\newcounter{subsub}[subsection]
\def\thesubsub{\thesubsection .\arabic{subsub}}
\def\subsub#1{\addtocounter{subsub}{1}\par\vspace{3mm}
\noindent{\bf \thesubsub ~ #1 }\par\vspace{2mm}}
\def\coker{\mathop{\rm coker}\nolimits}
\def\pr{\mathop{\rm pr}\nolimits}
\def\im{\mathop{\rm Im}\nolimits}
\def\hfl#1#2{\smash{\mathop{\hbox to 12mm{\rightarrowfill}}
\limits^{\scriptstyle#1}_{\scriptstyle#2}}}
\def\vfl#1#2{\llap{$\scriptstyle #1$}\big\downarrow
\big\uparrow
\rlap{$\scriptstyle #2$}}
\def\diagram#1{\def\normalbaselines{\baselineskip=0pt
\lineskip=10pt\lineskiplimit=1pt}   \matrix{#1}}
\def\limind{\mathop{\oalign{lim\cr
\hidewidth$\longrightarrow$\hidewidth\cr}}}

\long\def\InsertFig#1 #2 #3 #4\EndFig{
\hbox{\hskip #1 mm$\vbox to #2 mm{\vfil\includegraphics{#3}}#4$}}
\long\def\LabelTeX#1 #2 #3\ELTX{\rlap{\kern#1mm\raise#2mm\hbox{#3}}}

{\let\thefootnote\relax
\footnote{%\hskip3em 
\textbf{A.M.S.~classification :} 14J45, 14E30, 14E05, 32J18. 
}}

\vspace{-1cm}

\begin{center}Pr\'epublication de l'Institut Fourier n$^0$ 493 (2000) \\
{\em http://www-fourier.ujf-grenoble.fr/prepublications.html}
\end{center}

\bigskip

\bigskip

{\bf Abstract. } We prove finiteness of the number of 
smooth blow-downs on Fano manifolds
and boundedness 
results for the geometry of
non projective Fano-like manifolds.
Our proofs use properness of Hilbert schemes
and Mori theory.

\section*{Introduction}

In this Note, we say that a compact complex
manifold $X$ is a {\em Fano-like manifold}
if it becomes Fano after a finite sequence of blow-ups along
smooth connected centers, {\em i.e} 
if there exist a Fano manifold 
$\tilde{X}$ and a finite sequence of blow-ups along
smooth connected centers $\pi~: \tilde{X} \to X$. 
We say that a Fano-like manifold $X$ is {\em simple} 
if there exists a smooth submanifold $Y$ of $X$
($Y$ may not be connected) such that the blow-up
of $X$ along $Y$ is Fano.
If $Z$ is a projective manifold, we call 
{\em smooth blow-down of $Z$ (with an $s$-dimensional
center)} a map $\pi$ and a manifold $Z'$ such that 
$\pi~: Z\to Z'$ is 
the blow-up of $Z'$ along a 
smooth connected submanifold (of dimension $s$).
We say that a smooth blow-down of $Z$ is
projective (resp. non projective) if
$Z'$ is projective (resp. non projective).

\medskip

It is well-known that any Moishezon manifold becomes
projective after a finite sequence
of blow-ups along
smooth centers. Our aim is to bound the 
geometry of {\em Moishezon manifolds becoming Fano
after one blow-up along a smooth center}, i.e the geometry of
{\em simple non projective
Fano-like manifolds}.

\medskip

Our results in this direction are the following, 
the simple proof of Theorem~1 has been communicated to us 
by Daniel Huybrechts.

\medskip

{\bf Theorem~1.} {\em Let $Z$ be a Fano manifold of 
dimension $n$.
Then, there is only a finite number 
of smooth blow-downs of $Z$.}

\medskip

Let us recall here that the assumption $Z$ Fano
is essential~: there are projective smooth surfaces
with infinitely many $-1$ rational curves, hence
with infinitely many smooth blow-downs.

\medskip

Since there is only a finite number of deformation types 
of Fano manifolds of dimension $n$ (see 
\cite{KMM92} and also \cite{Deb97} for
a recent survey on Fano manifolds) and since smooth
blow-downs are stable under deformations \cite{Kod63}, 
we get the following
corollary (see section~1 for a detailed proof)~:

\medskip

{\bf Corollary~1.} {\em  There is only a finite 
number of deformation types
of simple Fano-like manifolds of dimension $n$.
}

\medskip

The next result is 
essentially due 
to Wi\'sniewski (\cite{Wis91}, prop.~(3.4) and (3.5)).
Before stating it, let us define  
$$ d_n = \max\{ (-K_Z)^n \, | \, Z \,\mbox{is a Fano manifold of dimension}
\, n  \}$$
and 
$$ \rho_n =  \max\{\rho (Z) := \rk (\Pic (Z)/\Pic ^0(Z)) 
\, | \, Z \,\mbox{is a Fano manifold of dimension} 
\, n  \}.$$

The number $\rho_n$ is well defined since there is only
a finite number 
of deformation types 
of Fano manifolds of dimension $n$ and
we refer to \cite{Deb97} for an explicit bound
for $d_n$.

\medskip

{\bf Theorem~2.} {\em Let $X$ be an $n$-dimensional simple 
non projective Fano-like manifold, $Y$ a smooth  
submanifold such that the blow-up $\pi~: \tilde{X}\to X$
of $X$ along $Y$ is Fano, and $E$ the exceptional
divisor of $\pi$. Then 
\begin{enumerate}
\item[(i)] if each component of $Y$ has Picard
number equal to one, 
then each component of $Y$ has ample
conormal bundle in $X$ and is Fano. Moreover
$\deg _{-K_{\tilde{X}}} (E) \leq (\rho_n -1) d_{n-1}$.
\item[(ii)] if $Y$ is a curve, then (each component of) $Y$
is a smooth rational curve
with normal bundle ${\mathcal O}_{\PP ^1}(-1)^{\oplus n-1}$. 
\end{enumerate}
}

\medskip

Finally, we prove here the following result~: 

\medskip

{\bf Theorem~3.} {\em Let $Z$ be a Fano manifold of 
dimension $n$ and index $r$. 
Suppose there is a non projective smooth blow-down
of $Z$ with an $s$-dimensional
center. Then $$r \leq (n-1)/2 \, \mbox{ and } \, s \geq r.$$
Moreover, 
\begin{enumerate}
\item[(i)] if $r > (n-1)/3$, then $s=n-1-r$~;
\item[(ii)] if $r < (n-1)/2$ and $s=r$, then $Y \simeq \PP ^r$.
\end{enumerate}
}

\medskip

Recall that the index of a Fano manifold $Z$
is the largest integer $m$ such that
$-K_Z = mL$ for $L$ in the Picard group of $Z$.

\medskip

{\bf Remarks.} 
\begin{enumerate}
\item[a)] For a Fano manifold $X$ of dimension $n$ and index $r$ with
second Betti number greater than or equal to $2$,
it is known that $\displaystyle{2r \leq n+2}$ \cite{Wi91},
with equality if and only if $X\simeq \PP ^{r-1} \times \PP ^{r-1}$.
\item[b)] Fano manifolds of even dimension (resp. odd dimension
$n$) and middle index (resp. index $(n+1)/2$) with $b_2 \geq 2$ 
have been intensively studied, see for example \cite{Wis93}.
Our Theorem~3 shows that there are no non projective smooth blow-down
of such a Fano manifold, without using any explicit classification.
\item[c)] The assumption that there is a {\em non projective} 
smooth blow-down
of $Z$ is essential in Theorem~3~: the Fano manifold
obtained by blowing-up $\PP ^{2r-1}$ along a $\PP ^{r-1}$
has index $r$.
\end{enumerate}
  
\section{Proof of Theorem~1 and Corollary~1. An example.}

\subsection{Proof of Theorem~1.}
{\em Thanks to D. Huybrechts for the following proof.}

Let $Z$ be a Fano
manifold and $\pi~: Z \to Z'$ a smooth blow-down 
of $Z$ with an $s$-dimensional
connected center.
Let $f$ be a line contained in a non trivial 
fiber of $\pi$. Then, the Hilbert polynomial
$P_{-K_Z}(m) := \chi (f,m(-K_Z)_{|f})$
is determined by $s$ and $n$ since $-K_Z \cdot f = n-s-1$
and $f$ is a smooth rational curve. 
Since $-K_Z$ is ample, the Hilbert scheme $\Hilb _{-K_Z}$
of
curves in $Z$ having $P_{-K_Z}$ as Hilbert polynomial
is a projective scheme, hence has a finite number
of irreducible components. Since each curve being 
in the component ${\mathcal H}$
of $\Hilb _{-K_Z}$ containing $f$ 
is contracted by $\pi$, there is only a finite
number of smooth blow-downs 
of $Z$ with an $s$-dimensional center.\finpreuve    

\subsection{Proof of Corollary~1.}
Let us first recall (\cite{Deb97} section 5.2) that there exists
an integer $\delta (n)$ such that every 
Fano $n$-fold can be realized as a smooth submanifold
of $\PP ^{2n+1}$ of degree at most $\delta (n)$.
Let us denote by $T$ a closed irreducible subvariety 
of the disjoint union of Chow varieties
of $n$-dimensional subvarieties of 
$\PP ^{2n+1}$ of degree at most $\delta (n)$,
and by $\pi : {\mathcal X}_T \to T$ the universal family.

{\em Step~1 : Stability of smooth blow-downs.}
Fix $t_0$ in the smooth locus $T_{smooth}$ 
of $T$ and suppose that $X_{t_0} := \pi ^{-1}(t_0)$ is a Fano $n$-fold  
and there exists a smooth blow-down of $X_{t_0}$
(denote by $E_{t_0}$ the exceptional 
divisor, $P$ its Hilbert polynomial with respect to 
${\mathcal O}_{\PP ^{2n+1}} (1)$). 
Let $S$ be the component of the Hilbert scheme  
of $(n-1)$-dimensional subschemes of 
$\PP ^{2n+1}$ with Hilbert polynomial $P$ and $u~: {\mathcal E}_S \to S$
the universal family.
Finally, let $I$ be the following subscheme of $T\times S$~:
$$ I = \{ (t,s) \, | \, u^{-1}(s) \subset X_t \}$$
and $p~: I \to T$ the proper algebraic map induced by the first
projection.
Thanks to the analytic stability of smooth blow-downs 
due to Kodaira (see \cite{Kod63}, Theorem~5),
the image $p(I)$ contains an analytic open neighbourhood
of $t_0$ hence it also contains a Zariski neighbourhood
of $t_0$. Moreover, since exceptional divisors are rigid,
the fiber $p^{-1}(t)$ is a single point for $t$ 
in a Zariski neighbourhood of $t_0$. Finally,
we get algebraic stability of smooth blow-downs 
(the $\PP ^r$-fibered structure of exceptional divisor
is also analytically stable - \cite{Kod63}, Theorem~4 -
hence algebraically stable by the same kind of argument).

\medskip

{\em Step~2 : Stratification of $T$ by the number
of smooth blow-downs}. 
For any integer $k\geq 0$, let us 
define
$$
\begin{array}{c} U_k(T) = \{ t \in T_{smooth} \, | \, X_t \,
\mbox{ is a Fano manifold and there exists 
at least}\\
 \mbox{  $k$ smooth blow-downs of } X_t \} ; 
\end{array}
$$ 
and $U_{-1}(T) = T_{smooth}$.
Thanks to Step~1, $U_k(T)$ is Zariski open in $T$, 
and thanks to Theorem~1, 
$$  \bigcap _{k\geq -1} U_k(T) = \emptyset.$$
Since $\{U_k(T)\}_{k \geq -1}$ 
is a decreasing sequence of Zariski open sets, 
by noetherian induction, we get that there
exists an integer $k$ such that $U_k(T)= \emptyset$
and we can thus define
$$ k(T) := \max \{ k \geq -1 \, | \, U_k(T) \neq \emptyset \}, \, \,
U(T):= U_{k(T)}(T).$$
Finally, we have proved that $U(T)$ is a non empty Zariski
open set of $T_{smooth}$ such that 
for every $t\in U(T)$, $Z_t$ is a Fano $n$-fold with
exactly $k(T)$ smooth blow-downs ($k(T)=-1$
means that for every $t \in T_{smooth}$, $X_t$ is
not a Fano manifold). 

Now let $T_0=T$, and $T_1$ be any closed
irreducible component of $T_0 \setminus U(T_0)$.
We get $U(T_1)$ as before and denote by $T_2$
any closed
irreducible component of $T_1 \setminus U(T_1)$, and so-on.
Again by noetherian induction, this process terminates after 
finitely many steps and we get a finite stratification
of $T$ such that each strata corresponds to
an algebraic family of Fano $n$-folds
with the same number of smooth blow-downs.  

\medskip

{\em Step~3 : Conclusion.} Since there is only a finite 
number of irreducible components in the Chow variety of Fano
$n$-folds, each being finitely stratified by Step~2, 
we get a finite number of deformation types
of simple Fano-like $n$-folds.\finpreuve

\medskip 

As it has been noticed by Kodaira, it is essential
to consider only {\em smooth} blow-downs. A $-2$ 
rational smooth curve on a surface is, in general, not stable
under deformations of the surface.  

\subsection{An example.}
Before going further, let us recall the following 
well known example.
Let $Z$ be the projective $3$-fold obtained by
blowing-up $\PP ^3$ along a smooth curve of type $(3,3)$
contained in a smooth quadric ${\mathcal Q}$
of $\PP ^3$. Let $\pi$ denotes the blow-up
$Z \to \PP^3$. Then $Z$ is a Fano manifold
of index one and there are at least three smooth blow-downs
of $Z$~: $\pi$, which is projective, and
two non projective smooth blow-downs consisting in contracting
the strict transform ${\mathcal Q}'$ of the quadric ${\mathcal Q}$
along one of its two rulings (the normal bundle
of ${\mathcal Q}'$ in $Z$ is ${\mathcal O}(-1,-1)$).

\medskip

{\bf Lemma~1.} {\em There are exactly three 
smooth blow-downs
of $Z$.}

\medskip

{\bf Proof~: }the Mori cone $\NE (Z)$ is a $2$-dimensional 
closed cone, one of its two extremal rays being generated
by the class of a line $f_{\pi}$ contained 
in a non trivial fiber of $\pi$,
the other one, denoted by $[R]$, by the class of one 
of the two rulings of  
${\mathcal Q}'$ (the two rulings are numerically
equivalent, the corresponding extremal contraction
consists in contracting ${\mathcal Q}'$ to a singular point
in a projective variety, hence is not a smooth blow-down).
If $E$ is the exceptional divisor of $\pi$, we have
$$ E\cdot [f_{\pi}] = -1, \, E\cdot [R] = 3, 
\, {\mathcal Q}'\cdot [f_{\pi}] =1,
\, {\mathcal Q}'\cdot [R] = -1. $$
Now suppose there exists a smooth blow-down $\tau$ 
of $Z$ with a $1$-dimensional center, which is not one 
of the three previously described. Let  
$L$ be a line contained in a non trivial fiber of $\tau$,
then since $-K_Z\cdot [L] = 1$, 
we have 
$[L] = a [f_{\pi}] + b[R]$ for some strictly positive numbers
such that $a+b =1$.
Since we have moreover 
$$ {\mathcal Q}'\cdot  [L] = a-b = 2a -1 \in \ZZ
\,\, \mbox{and} \,\, E \cdot [L] = 3b-a = 3 - 4a,$$ 
we get $a=b=1/2$. 
Therefore ${\mathcal Q}'\cdot  [L] =0$ hence $L$ is disjoint
from ${\mathcal Q}'$ (it can not be contained in ${\mathcal Q}'$
since ${\mathcal Q}'_{|{\mathcal Q}'} = {\mathcal O}(-1,-1)$).
It implies that there are two smooth blow-downs
of $Z$ with disjoint exceptional divisors, which is 
impossible since $\rho (Z) = 2$.

Finally, if there is a smooth blow-down $\tau~: Z \to Z'$ 
of $Z$ with a $0$-dimensional center, then $Z'$ is 
projective and $\tau$ is a Mori extremal contraction,
which is again impossible since we already met the two
Mori extremal contractions on $Z$. \finpreuve 

\section{Non projective smooth blow-downs on a center with Picard number
$1$.
Proof of Theorem~2.}

The proof of Theorem~2 we will give
is close to 
Wi\'sniewski's one but we give two intermediate results
of independant interest.  

\subsection{On the normal bundle of the center.}
Let us recall that a smooth submanifold $A$ in a complex
manifold $W$ is contractible to a point (i.e.
there exists a complex space $W'$ and a map $\mu~: W \to W'$  
which is an isomorphism outside $A$ and such
that $\mu (A)$ is a point) if and only if
$N_{A/W}^*$ is ample (Grauert's criterion \cite{Gra62}). 

The following proposition was proved by Campana \cite{Cam89}
in the case where $Y$ is a curve and $\dim (X) =3$.

\medskip

{\bf Proposition~1.} {\em Let $X$ be a
non projective manifold, $Y$ a smooth
submanifold of $X$ such that the blow-up $\pi~: \tilde{X}\to X$
of $X$ along $Y$ is projective.
Then, for each connected component $Y'$ of
$Y$
with $\rho (Y') =1$,
the conormal bundle $N_{Y'/X}^*$ is ample.
}

\medskip

Before the proof, let us remark that $Y$ is 
projective since the exceptional divisor of $\pi$ is.

\medskip

{\bf Proof of Proposition~1~: }(following Campana) we can suppose that 
$Y$ is connected. Let $E$ be
the exceptional divisor of $\pi$
and $f$ a line contained in a non trivial 
fiber of $\pi$. Since $E \cdot f = -1$, there is
an extremal ray $R$ of the Mori cone $\overline{\NE} (\tilde X)$
such that $E\cdot R < 0$. 
Since $E\cdot R < 0$, $R$ defines 
an extremal ray of the Mori cone 
$\overline{\NE} (E)$ which we still denote
by $R$ (even if $\overline{\NE} (E)$ 
is not a subcone of $\overline{\NE} (\tilde X)$
in general !).
Since $\rho (Y) =1$, we have $\rho (E) =2$,
hence $\overline{\NE} (E)$ is a $2$-dimensional 
closed cone, one of its two extremal rays being generated
by $f$. Then~:
\begin{enumerate}
\item[-] either $R$ is not generated by $f$
and $E_{|E}$ is strictly negative on $\overline{\NE} (E)\setminus \{ 0 \}$.
In that case, $-E_{|E} = {\mathcal O }_E (1)$ 
is ample by Kleiman's criterion,
which means that $N_{Y/X}^*$ is ample.
\item[-] or, $R$ is generated by $f$. In that case, the
Mori contraction $\varphi _R ~: \tilde X \to Z$
fac\-to\-ri\-ze through $\pi$~:

\centerline{
\xymatrix{ \tilde{X} \ar[d]_{\pi} \ar[r]^{\varphi _{R}}
  & Z\\
X \ar[ur]_{\psi} & 
}
}

\noindent where $\psi~:X\to Z$ is an isomorphism
outside $Y$. Since 
the variety $Z$ is projective
and $X$ is not,
$\psi$ is not an isomorphism 
and since $\rho (Y) =1$, $Y$ is contracted to a point by $\psi$,
hence $N_{Y/X}^*$ is ample by Grauert's criterion.\finpreuve
\end{enumerate}

\medskip

Let us prove the following consequence of Proposition~1:

\medskip

{\bf Proposition~2.} {\em Let $X$ be a
non projective manifold, $Y$ a smooth
submanifold of $X$
such that the blow-up $\pi~: \tilde{X}\to X$
of $X$ along $Y$ is projective with $-K_{\tilde X}$
numerically effective (nef).
Then, each connected component $Y'$ of
$Y$
with $\rho (Y') =1$
is a Fano manifold.
}

\medskip

{\bf Proof~: }we can suppose that 
$Y$ is connected. Let $E$ be
the exceptional divisor of $\pi$. 
Since $- E_{|E}$
is ample by Proposition~1, the adjunction formula
$-K_E = -K_{\tilde X | E} - E_{|E}$
shows that $-K_E$ is ample, hence $E$ is Fano.
By a result of Szurek and Wi\'sniewski \cite{SzW90}, $Y$
is itself Fano.\finpreuve

\subsection{Proof of Theorem~2.}

For the first assertion, we only have
to prove that 
$$\deg _{-K_{\tilde{X}}} (E) \leq (\rho_n -1) d_{n-1}.$$
Let $Y'$ be a connected component of $Y$ and
$E' = \pi^{-1}(Y')$. 
Then, since $-E_{|E'}$ is ample~:
$$ \deg_{-K_{\tilde X}}(E') =
(-K_{\tilde {X}|E'})^{n-1} = (-K_{E'}+ E_{|E'})^{n-1}
\leq  (-K_{E'})^{n-1} \leq d_{n-1}.$$
Now, if $m$ is the number of connected components of $Y$,
then $$\rho(\tilde X) = m +\rho (X) \geq m+1.$$
Putting all together, we get
$$\deg _{-K_{\tilde{X}}} (E) \leq  (\rho_n -1) d_{n-1},$$
which ends the proof of the first point.

We refer to \cite{Wis91} prop.~(3.5) for the second point.\finpreuve

\section{On the dimension of the center of non projective smooth
blow-downs.
Proof of Theorem~3.}

Theorem~3 is a by-product of the more precise
following statement and of Proposition~3 below~:

\medskip

{\bf Theorem~4. } {\em Let $Z$ be a Fano manifold of
dimension $n$ and index $r$,
$\pi~: Z \to Z'$ be a non projective
smooth blow-down of $Z$, $Y \subset Z'$
the center of $\pi$. Let $f$ be a line contained in a non trivial
fiber of $\pi$, then 
\begin{enumerate}
\item [(i)] if $f$ generates an extremal ray of $\NE (Z)$,
then $\dim (Y) \geq (n-1)/2$.
\item [(ii)] if $f$ does not generate an extremal ray of $\NE (Z)$,
then $\dim (Y) \geq r$. Moreover, 
if $\dim (Y) = r$, then  
$Y$ is isomorphic to $\PP ^r$. 
\end{enumerate}
In both cases (i) and (ii), $Y$ contains 
a rational curve. 
}

\medskip

The proof relies on Wi\'sniewski's inequality (see \cite{Wis91}
and \cite{AnW95}), which we recall
now for the reader's convenience~: let $\varphi~: X \to Y$
be a Fano-Mori contraction (i.e $-K_X$ is $\varphi$-ample)
on a projective manifold $X$, $\Exc (\varphi)$ its exceptional
locus and
$$ l({\varphi}) := \min \{ -K_X\cdot C \,;\, C \, \,
\mbox{rational curve contained in}
\Exc (\varphi) \}$$
its length, then for every non trivial fiber $F$~:
$$ \dim \Exc (\varphi) + \dim (F) \geq \dim (X) -1 + l({\varphi}).$$

\medskip 

{\bf Proof of Theorem~4.} The method of proof is taken from
Andreatta's recent paper \cite{And99} (see also \cite{Bon96}). 

{\em First case~:} suppose that a line $f$ contained in a non trivial
fiber of $\pi$ generates an extremal ray $R$ of $\NE (Z)$.
Then the 
Mori contraction 
$\varphi _R ~: Z \to W$
fac\-to\-ri\-zes through $\pi$~:
\centerline{
\xymatrix{ Z \ar[d]_{\pi} \ar[r]^{\varphi _{R}}
  & W\\
Z' \ar[ur]_{\psi} & 
}
} 
\noindent where $\psi$ is an isomorphism outside
$Y$. In particular, the exceptional locus $E$
of $\pi$ is equal to the exceptional locus
of the extremal contraction $\varphi _{R}$. 

Let us now denote by $\psi _Y$ the restriction of
$\psi$ to $Y$, $s=\dim (Y)$, $\pi _E$ and 
$\varphi _{R,E}$ the restriction of $\pi$
and $\varphi_{R}$ to $E$. Since $Z'$ is not projective, 
$\psi _Y$ is not a finite map.
Since $\varphi _{R}$ is birational, 
$W$ is $\QQ$-Gorenstein, hence 
$K_W$ is $\QQ$-Cartier and $K_{Z'}=\psi ^* K_W$. Therefore,
$K_{Z'}$
is $\psi$-trivial, 
hence 
$K_Y + \det N_{Y/Z'}^*$ is $\psi _Y$-trivial.
Moreover,
${\mathcal O}_E (1)= -E_{|E}$ 
is $\varphi _{R,E}$-ample by Kleiman's criterion,
hence $N_{Y/Z'}^*$ is $\psi_Y$-ample.
Finally, $\psi_Y$ is a Fano-Mori contraction,
of length greater or equal to $n-s =\rk ( N_{Y/Z'}^*)$.
Together with Wi\'sniewski's inequality 
applied on $Y$, we get that for every non trivial fiber $F$ 
of $\psi _Y$
$$ 2s \geq \dim (F) + \dim \Exc (\psi _Y) \geq 
n-s + s -1$$
hence $2s \geq n-1$.
Moreover, $\Exc (\psi _Y)$ is covered by rational curves,
hence $Y$ contains a rational curve. 

{\em Second case~:} suppose that a line $f$ contained in a non trivial
fiber of $\pi$ does not generate an extremal ray $R$ of $\NE (Z)$.
In that case, since $E \cdot f =-1$, there is an extremal ray
$R$ of $\NE (Z)$ such that $E\cdot R < 0$. In particular,
the exceptional locus $\Exc (R)$ of the extremal contraction
$\varphi _{R}$ is contained in $E$, and
since $f$ is not on $R$, we get for
any fiber $F$ of $\varphi _{R}$~:
$$\dim ( F) \leq s = \dim (Y).$$

By the adjunction formula, $-K_E = - K_{Z|E} - E_{|E}$,
the length $l_E(R)$ of $R$ as {\em an extremal ray of $E$}
satisfies 
$$ l_E(R) \geq r +1,$$
where $r$ is the index of $Z$.
Together with Wi\'sniewski's inequality
applied on $E$, we get~:
$$ r+1 + (n-1)-1 \leq s + \dim (\Exc (R)) \leq s +n-1.$$ 
Finally, we get $r\leq s$, and since the fibers
of $\varphi _{R}$ are covered by rational curves,
there is a rational curve in $Y$. 
Suppose now (up to the end) that $r=s$. 
Then $E$ is the exceptional locus of the 
Mori extremal contraction $\varphi _{R}$.
Moreover, $K_Z + r(-E)$ is 
a good supporting divisor for $\varphi _{R}$, and since 
every non trivial fiber of $\varphi _{R}$ has dimension $r$,
$\varphi _{R}$ is a smooth projective blow-down.
In particular, the restriction of $\pi$ to a non 
trivial fiber $F\simeq \PP ^r $ induces a 
finite surjective map $\pi~: F\simeq \PP ^r \to Y$
hence $Y\simeq \PP ^r$ by a result of Lazarsfeld \cite{Laz83}.

This ends the proof of Theorem~4. 
\finpreuve

\medskip

The proof of Theorem~4 does not use the hypothesis
$Z$ Fano in the first case. We therefore have the
following~:

\medskip

{\bf Corollary~2.} {\em Let $Z$ be a projective manifold of
dimension $n$,  
$\pi~: Z \to Z'$ be a non projective
smooth blow-down of $Z$, $Y \subset Z'$
the center of $\pi$. Let $f$ be a line contained in a non trivial
fiber of $\pi$ 
and suppose $f$ generates an extremal ray of $\overline{\NE} (Z)$.
Then $\dim (Y) \geq (n-1)/2$. Moreover, if $\dim (Y)= (n-1)/2$,
then $Y$ is contractible on a point.
}

\medskip 

We finish this section by the following easy proposition, 
which combined with 
Theorem~4 implies Theorem~3 of the Introduction~:

\medskip

{\bf Proposition~3.} {\em Let $Z$ be a Fano manifold of
dimension $n$ and index $r$,
$\pi~: Z \to Z'$ be a smooth blow-down of $Z$, $Y \subset Z'$
the center of $\pi$. Then $n-1-\dim (Y)$ is a multiple of $r$.
}

\medskip

{\bf Proof.} Write $$-K_Z = rL \, \, \mbox{ and } \,\, 
-K_Z = -\pi ^* K_{Z'}-(n-1-\dim (Y))E$$
where $E$ is the exceptional divisor of $\pi$.
Let $f$ be a line contained in a fiber of $\pi$.
Then $rL\cdot f = n-1- \dim (Y)$, which ends the proof.\finpreuve  

\medskip 

{\bf Proof of Theorem~3.}
Let $Z$ be a Fano manifold of 
dimension $n$ and index $r$ and  
suppose there is a non projective smooth blow-down
of $Z$ with an $s$-dimensional
center. By Proposition~3, there is a strictly positive integer
$k$ such that $n-1- kr = s$.
By Theorem~4, either $n-1-kr \geq (n-1)/2$
or $n-1-kr \geq r$. In both cases, it implies that
$r\leq (n-1)/2$ and therefore $s\geq r$.
If $r > (n-1)/3$, since $n-1 \geq (k+1)r > (k+1)(n-1)/3$, we get
$k=1$ and $s=n-1-r$.\finpreuve

\section{Rational curves on simple Moishezon manifolds.}

\medskip

The arguments of the previous section
can be used to deal
with the following well-known question~:
{\em does every non projective Moishezon manifold
contain a rational curve ?}
The answer is positive in dimension three (it is 
due to Peternell \cite{Pet86}, see also \cite{CKM88} p. 49
for a proof using the completion
of Mori's program in dimension three).   

\medskip

{\bf Proposition~4.} {\em Let $Z$ be a projective manifold,
$\pi~: Z \to Z'$ be a non projective smooth blow-down of $Z$.
Then $Z'$ contains a rational curve.
}

\medskip 

{\bf Proof.} With the notations of the previous section, it
is clear in the first case where a line $f$ contained in a non trivial
fiber of $\pi$ generates an extremal ray $R$ of $\overline{\NE} (Z)$
(in that case, the center of $\pi$ contains a rational curve).
In the second case, since $f$ is not extremal 
and $K_Z$ is not nef, there is a Mori contraction $\varphi$ on $Z$
such that any rational curve contained in a fiber of $\varphi$
is mapped by $\pi$ to a non constant rational curve in $Z'$.\finpreuve

-----------

{\em 
\noindent L.B.~: Institut Fourier, UMR 5582, 
Universit\'e de Grenoble 1, BP 74. 
38402 Saint Martin d'H\`eres. FRANCE\\
\noindent e-mail : bonavero@ujf-grenoble.fr
}

\medskip

{\em
\noindent S.T.~: Department of Mathematics,
Graduate School of Science, Osaka University.
Toyonaka, Osaka, 560-0043 JAPAN.\\
\noindent e-mail : taka@math.sci.osaka-u.ac.jp
}

\end{document}